\documentstyle{article}

\input amssym.def
\input amssym

\newcommand{\F}{{\cal F}} 
\newcommand{\ar}{\rightarrow}
  
\newcommand{\Z}{{\Bbb Z}} 
\newcommand{\R}{{\Bbb R}}

\newcommand{\Cinf}{C^{\infty}}
\newtheorem{thm}{Theorem}[section]

\newtheorem{cor}[thm]{Corollary}

\newtheorem{prop}[thm]{Proposition}
\newtheorem{example}[thm]{Example}

\newcommand{\be}{\begin{equation}}
\newcommand{\ee}{\end{equation}}
\newcommand{\bea}{\begin{eqnarray}}
\newcommand{\eea}{\end{eqnarray}}
\newcommand{\ba}{\begin{eqnarray*}}
\newcommand{\ea}{\end{eqnarray*}}
\newcommand{\bt}{\begin{thm}}
\newcommand{\et}{\end{thm}}
\newcommand{\bc}{\begin{cor}}
\newcommand{\ec}{\end{cor}}
\newcommand{\bp}{\begin{prop}}
\newcommand{\ep}{\end{prop}}
\newcommand{\cC}{{\cal C}}
\newcommand{\cH}{{\cal H}}
\newcommand{\g}{{\frak g}}

\begin{document}
\title{On the first secondary invariant of Molino's central 
sheaf\thanks{Supported in part by Spanish DGICYT grant
PB90-0765}} 
\author{Jes\'us A. \'ALVAREZ L\'OPEZ\\[6pt]
\small Universidade de Santiago de Compostela,\\ 
\small Facultade de Matem\'aticas, Campus Sur,\\
15782 Santiago de Compostela,  Spain\\
\small {\bf E-mail} jesus.alvarez@usc.es } 
\date{} 
\maketitle

\begin{abstract} For a Riemannian foliation on a closed manifold,
the first secondary invariant of Molino's central sheaf is
an obstruction to tautness. Another obstruction is the
class defined by the basic component of the mean curvature with respect to some
metric. Both obstructions are proved to be the same up to a constant, and other
geometric properties are also proved to be equivalent to tautness. \end{abstract}

\section{Introduction and main results}\label{sec:intro}
Let $\F$ be a Riemannian foliation on a closed manifold $M$ \cite{Reinhart},
$\Omega^\cdot(M/\F)$ its basic complex, and $H^\cdot(M/\F)$ its basic cohomology
\cite{ElKacimiHector,ElKacimiHectorSergiescu,KamberTond2}. There is a locally trivial sheaf
$\cC=\cC(\F)$ of Lie algebras of germs of transverse Killing fields whose `transverse orbits'
are the leaf closures  \cite{Molino1,Molino2}. It is called the central sheaf of $\F$. The
typical fiber of $\cC$ is the opposite of the structural Lie algebra $\g$ of $\F$. The sheaf
$\cC$ canonically defines a vector bundle $C=C(\F)$ over $M$ with a flat connection. The
corresponding multiplicative homomorphism $\Delta_\ast:H^\cdot({\rm gl(q)},{\rm O}(q))\ar
H^\cdot(M)$ \cite{KamberTond1}, $q={\rm codim}\,\F$, can be given as a composite of a
homomorphism $\Delta_\ast=\Delta(\F)_\ast:H^\cdot({\rm gl}(q),{\rm O}(q))\ar H^\cdot(M/\F)$
and the canonical homomorphism $H^\cdot(M/\F)\ar  H^\cdot(M)$. We get basic secondary
invariants $\Delta_\ast(y_i)=\Delta(\F)_\ast(y_i)\in H^{2i-1}(M/\F)$,
$i=1,\ldots,2[(m+1)/2]-1$, $m=\dim\g$. The basic class $\Delta_\ast(y_1)$ will be studied in
this paper. It would be also interesting to study the geometric information contained in the
$\Delta_\ast(y_i)$ for $i>1$. 

 It was pointed 
out in \cite{MolinoSergiescu} that, for Riemannian flows, $\Delta_\ast(y_1)$ is the
obstruction to tautness; i.e. the obstruction to the existence of a metric on $M$
such that the leaves are minimal submanifolds. This property also holds for $\F$ of
arbitrary dimension: Suppose $\F$ is transversely orientable for simplicity, then $\F$
is taut if and only if $H^q(M/\F)\neq 0$ \cite{Masa2}, which is equivalent to the triviality
of the sheaf $\bigwedge^m\cC$ \cite{Sergiescu}, and this in turn is equivalent to
$\Delta_\ast(y_1)=0$. There is another obstruction to tautness: For any bundle-like metric,
the basic component of the mean curvature form of the leaves is closed and defines a class
$\xi=\xi(\F)\in H^1(M/\F)$, which depends only on $\F$ and vanishes if and only if $\F$ is
taut \cite{Alv7}. We close this circle of ideas by proving directly that both obstructions
are the same up to a constant:
 
\bt\label{thm:xi} With the above notation, $\xi=-2\pi\,\Delta_\ast(y_1)$.\et

We also prove other relations between $\Delta_\ast(y_1)$ and geometric properties of
$\F$.  Consider the filtration of $\Omega^\cdot(M/\F)$ given by the differential ideals
$F^k\Omega^\cdot(M/\F)$, where an $\alpha\in\Omega^r(M/\F)$ is in  $F^k\Omega^r(M/\F)$ if
$ i_X\alpha=0$ for $X=X_1\wedge\ldots\wedge X_{r-k+1}$ with the vector fields $X_j$
tangent to the leaf closures. The corresponding spectral sequence $(E_i,d_i)$ converges to
$H^\cdot(M/\F)$ (cf. \cite[\S 2]{ElKacimiNicolau}). If $M/\bar{\F}$ denotes the space of
leaf closures of $\F$, there is a canonical isomorphism  $E_2^{\cdot,0}\cong
H^\cdot(M/\bar{\F})$. So there is a canonical injection  $H^1(M/\bar{\F})\ar
H^1(M/\F)$. Let $F^kH^\cdot(M/\F)$ be the induced filtration of $H^\cdot(M/\F)$. The element
defined by $\Delta_\ast(y_1)$ in $H^1(M/\F)/F^1H^1(M/\F)\equiv E_\infty^{0,1}$ will be
denoted by $\bar{\Delta}_\ast(y_1)$. Thus $\bar{\Delta}_\ast(y_1)=0$ if  and only if
$\Delta_\ast(y_1)\in F^1H^1(M/\F)\equiv E_\infty^{1,0}\cong H^1(M/\bar{\F})$.
 
\bt With the above notation, $\bar{\Delta}_\ast(y_1)=0$ if and only if $\g$ is unimodular.\et

\bt With the above notation, suppose $\F$ admits a transverse parallelism. Let $\cH$ be any
representative of the holonomy pseudogroup of $\F$ on some manifold $T$. Then: 
\begin{itemize}
\item[(i)] If $\g$ is unimodular, $\Delta_\ast(y_1)=0$ if and only if the $\cH$-orbit closures
are minimal submanifolds for some $\cH$-invariant metric on $T$.
\item[(ii)] If $\g$ is not unimodular, the $\cH$-orbit closures are minimal submanifolds for
some $\cH$-invariant metric on $T$. 
\end{itemize}
Thus $\F$ is taut if and only if $\g$ is unimodular and the $\cH$-orbit closures are minimal
submanifolds for some $\cH$-invariant metric. \et

If $\F$ does not admit any transverse parallelism, a similar result can be stated by
considering the horizontal lifting $\hat{\F}$ to the principal bundle of transverse
orthonormal frames for some fixed transverse Riemannian structure \cite{Molino1,Molino2}. 
In particular we have the
following. 

\bc With the above notation, let $\hat{\cH}$ be the holonomy pseudogroup of $\hat{\F}$. Then 
$\F$ is taut if and only if $\g$ is unimodular and the $\hat{\cH}$-orbit closures are minimal
submanifolds for some $\hat{\cH}$-invariant metric. \ec

For a bundle-like metric on $M$, let $\kappa$ be the mean curvature form of the leaves, and
$\kappa_b$ its basic component. 

\bt With the above notation, the following holds:
\begin{itemize}
\item[(i)] If $\g$ is unimodular then, for any bundle-like metric, $\kappa_b$ vanishes on
vectors tangent to the leaf closures.  
\item[(ii)] If $\g$ is not unimodular, there is a
bundle-like metric such that $\kappa_b$ vanishes on vectors orthogonal to the leaf
closures.
\end{itemize}\et

This theorem can be sharpened by the recent result of D.~Dom\'{\i}nguez \cite{Demetrio},  
showing the existence of a bundle-like metric on $M$ with basic mean curvature of the
leaves. Indeed any representative of $\xi$ can be realized as the mean curvature for some
bundle-like metric. For Lie foliations with dense leaves, the result is very explicit:

\bc If $\F$ is a Lie $\g$-foliation with dense leaves then, for any bundle-like metric, the
basic component of the mean curvature corresponds to the trace of the adjoint representation
by the canonical identity $\Omega^\cdot(M/\F)\equiv \bigwedge^\cdot\g^\ast$. Moreover such
form can be always realized as the mean curvature for some bundle-like metric.\ec

These results depend only on the holonomy pseudogroup of $\F$. Thus, with slightly 
more generality, we shall consider a complete pseudogroup $\cH$ of local isometries of a
Riemannian manifold $T$ \cite{Haefliger85,Haefliger88}.

\section{A remark on the first secondary characteristic class of a flat vector
bundle}\label{sec:remark} Let $(M,g)$ be a Riemannian manifold, $\rho:V\ar M$ a rank
$m$ vector bundle with a flat connection $\nabla$, $P$ the ${\rm Gl}(m)$-principal bundle of
frames of $V$, and $\omega$ the connection form defined by $\nabla$. The induced connection
on $\bigwedge V$ will be also denoted by $\nabla$. For any given ${\rm O}(m)$-reduction of $P$
defined by a section $s:M\ar P/{\rm O}(m)$, we have the multiplicative homomorphism
$\Delta_\ast:H({\rm gl}(m),{\rm O}(m))\ar H^\cdot(M)$ \cite[Theorem~4.43]{KamberTond1},
yielding secondary characteristic invariants $\Delta_\ast(y_i)\in H^{2i-1}(M)$ for
$i=1,\ldots,2[(m+1)/2]-1$ \cite[Theorem~6.33]{KamberTond1}. A representative of
$\Delta_\ast(y_1)$ is the form $\Delta(y_1)=\frac{1}{2\pi}\,s^\ast({\rm
trace}\,\omega)\in\Omega^1(M)$ \cite[Proof of Proposition~6.34]{KamberTond1}.

Recall that any smooth section $X$ of $\bigwedge V$ canonically defines a smooth section
$\hat{X}$ of $\rho^\ast\bigwedge V\equiv\bigwedge\rho^\ast V$. Identifying $\rho^\ast V$
with the vertical bundle of $\rho$ in the canonical way, we can consider such $\hat{X}$ as a
smooth section of $\bigwedge{\rm T}V$ over $V$. Moreover, if $\tilde{Z}$ is the horizontal
lifting of any vector field $Z$ on $M$, we get
\be\label{eqn:hatnabla} \widehat{\nabla_ZX}=\theta_{\tilde{Z}}\hat{X}\;,\ee
where $\theta_{\tilde{Z}}$ denotes Lie derivative with respect to $\tilde{Z}$. This can be
seen as follows. The parallel transport along the integral curves of $Z$ is given by the
integral curves of $\tilde{Z}$. In particular, the restriction of the flow of $\tilde{Z}$
between two fibers of $\rho$ is linear, and thus can be canonically identified with its
derivative at each point. Therefore Eq.~(\ref{eqn:hatnabla}) follows from the usual expression
of covariant derivative in terms of parallel transport and the usual expression of Lie
derivative in terms of the flow of vector fields.

Consider the Riemannian structure on $V$ defined by the ${\rm O}(m)$-reduction of $P$. 
We get an induced Riemannian structure on the
vertical bundle of $\rho$ by identifying it with $\rho^\ast V$. Let $\hat{g}$ be the
Riemannian metric on $V$ defined as the orthogonal sum of the lift of $g$ to the horizontal
bundle and the Riemannian structure  on the vertical bundle. The $\hat{g}$-mean curvature
form of the fibers of  $\rho$ will denoted by $\kappa_V$. The induced metric on
$\bigwedge{\rm T}V$ will be also denoted by $\hat{g}$. 

\bp\label{prop:kappa} $\kappa_V=2\pi\,\rho^\ast\Delta(y_1)$.\ep 

{\em Proof}.  We can clearly
assume $V$ is an oriented vector bundle. So $P$ has a ${\rm Gl}^+(m)$-reduction $P^+$.
Consider the homomorphism $\det : {\rm Gl}^+(m)\ar\R^+$, and the corresponding bundle map
$P^+\ar\bar{P}=P^+\times_{{\rm Gl}^+(m)}\R^+$. Then $\nabla$ defines a
flat connection on $\bar{P}$, and let $\bar{\omega}$ be its connection form.

The section $s$ defines a section $\bar{s}$ of $\bar{P}\ar M$ because the
composite 
$${\rm SO}(m)\hookrightarrow{\rm Gl}^+(m)\stackrel{\rm det}{\longrightarrow}\R^+$$ 
is trivial. By
functoriality of the construction of the characteristic homomorphism under homomorphisms of
structural groups \cite[Theorem~4.43~(iii)]{KamberTond1}, we have
\be\label{eqn:somega} \bar{s}^\ast\bar{\omega}=2\pi\,\Delta(y_1)\;.\ee 
(See the proof of Proposition~6.34 in \cite{KamberTond1}.)

Since $V$ is an oriented bundle, there is a non-vanishing section $X\in\Cinf(\bigwedge^mV)$,
with a corresponding section 
$\hat{X}\in\Cinf(\rho^\ast\bigwedge^mV)\equiv\Cinf(\bigwedge^m\rho^\ast V)$. By identifying 
$\rho^\ast V$ with the vertical bundle, if $X$ is
unitary, then $\chi=\hat{g}(\hat{X},\cdot)$ is the characteristic form of the fibers of $\rho$
\cite{Rummler}; i.e. $\chi(U)=\hat{g}(\hat{X},U)$ for any $U\in\Cinf(\bigwedge{\rm T}V)$. Thus
\be\label{eqn:Rummler} \theta_Y\hat{X}=\kappa_V(Y)\,\hat{X}\ee for any horizontal
$\rho$-projectable vector field $Y$ on $V$. Indeed $\theta_Y\hat{X}=f\,\hat{X}$ for some
function $f$ on $V$ because the flow of $Y$ maps fibers of $\rho$ to fibers of $\rho$, and
Rummler's mean curvature  formula implies \cite{Rummler} $$0=\theta_Y(\chi(\hat{X}))=
(\theta_Y\chi)(\hat{X})+\chi(\theta_Y(\hat{X}))=
-\kappa_V(Y)+f\;.$$

On the other hand, 
$\bar{P}$ can be canonically identified with the principal bundle of oriented frames of
the line bundle $\bigwedge^mV$. Thus 
\be\label{eqn:nabla}\nabla_ZX=(\bar{s}^\ast\bar{\omega})(Z)\,X\;.\ee 
for any vector field $Z$ on $M$, where $\tilde{Z}$ is its horizontal lifting.

Therefore
\bea
\kappa_V(\tilde{Z})\,\hat{X}&=&\theta_{\tilde{Z}}\hat{X}\quad\mbox{by
Eq.~(\ref{eqn:Rummler})}\nonumber\\
&=&\widehat{\nabla_ZX}\quad\mbox{by Eq.~(\ref{eqn:hatnabla})}\nonumber\\
&=&(\bar{s}^\ast\bar{\omega})(Z)\,\hat{X}
\quad\mbox{by Eq.~(\ref{eqn:nabla})}\nonumber\\
&=&2\pi\,\rho^\ast\Delta(y_1)(\tilde{Z})\,\hat{X}\quad\mbox{by Eq.~(\ref{eqn:somega})}\;.\nonumber
\eea
The result now follows because $\kappa_V$ and $\rho^\ast\Delta(y_1)$ vanish on
vertical vectors. $\Box$

\section{Preliminaries on complete pseudogroups of local
isometries}\label{sec:preliminaries} 
Let $\cH$ be a complete pseudogroup of local 
isometries of a Riemannian manifold $(T,g)$, $T/\cH$ the space of
$\cH$-orbits, and $\bar{\cH}$ the closure of $\cH$ \cite{Haefliger85}. Thus
$T/\bar{\cH}$ is the space of $\cH$-orbit closures.  

If $\cH$ preserves a parallelism on $T$, then we have the following description due to
E.~Salem \cite{Salem}. The space $T/\bar{\cH}$ is a manifold and the canonical
projection $\pi_b:T\ar T/\bar{\cH}$ is a submersion. Moreover, for some Lie group $G$
and some dense subgroup $\Lambda\subset G$, every point in $T/\bar{\cH}$ has a
neighborhood $U$ so that the restriction of $\cH$ to $\pi_b^{-1}(U)$ is equivalent to the
pseudogroup generated by the action of $\Lambda$ on $G\times U$, acting by left
multiplication on $G$ and trivially on $U$. Furthermore $\pi_b$ corresponds to the canonical
second projection of $G\times U$ onto $U$ by this equivalence. The Lie algebra $\g$ of $G$
is called the structural Lie algebra of $\cH$, and $\pi_b$ its basic projection. 

For arbitrary $\cH$, it is standard to consider the ${\rm O}(n)$-principal bundle
$\pi:\hat{T}\ar T$ of orthonormal frames on $T$ with the Levi-Civita connection, where 
$n=\dim T$, and the complete pseudogroup $\hat{\cH}$ canonically defined by $\cH$ on
$\hat{T}$. The canonical parallelisms on $\hat{T}$ are $\hat{\cH}$-invariant, thus Salem's
description holds for $\hat{\cH}$. The structural Lie algebra of $\hat{\cH}$ is also called
the structural Lie algebra of $\cH$. (There is no ambiguity when $\cH$ preserves a
parallelism.) The ${\rm O}(n)$-action on $\hat{T}$ preserves $\hat{\cH}$, and thus there is
an induced ${\rm O}(n)$-action on the manifold $W$ of $\hat{\cH}$-orbit closures so that the
basic projection $\pi_b$ is ${\rm O}(n)$-equivariant, yielding a canonical identity
\be\label{eqn:W} T/\bar{\cH}\equiv W/{\rm O}(n)\;.\ee

The complex of 
$\cH$-invariant differential forms will be denoted by
$\Omega^\cdot_\cH=\Omega^\cdot(T)_\cH$, and its cohomology by $H^\cdot(T)_\cH$. 
We shall also use the notation $\Omega^\cdot_{\cH, i =0}$ for the space of $\cH$-invariant
forms which vanish on vector fields tangent to the $\cH$-orbit closures. Similarly, let 
$\Omega(W)_{{\rm O}(n), i =0}$ be the complex of ${\rm O}(n)$-invariant differential forms
on $W$ which vanish on vector fields tangent to the ${\rm O}(n)$-orbits. 

Define a filtration
of $\Omega^\cdot_\cH$ by differential ideals $F^k\Omega^\cdot_\cH$ where an $\cH$-invariant
$r$-form $\alpha$ is in $F^k\Omega^r_\cH$ if $ i_X\alpha=0$ for $X=X_1\wedge\ldots\wedge
X_{r-k+1}$ with the vector fields $X_j$ tangent to the orbit closures. The corresponding
spectral sequence $(E_i,d_i)=(E_i(\cH),d_i)$ converges to $H^\cdot(T)_\cH$. We have
$E_0^{u,\cdot}=F^u\Omega^\cdot_\cH/F^{u+1}\Omega^\cdot_\cH$
with the differential map induced by the de~Rham derivative, and thus
$$E_1^{u,v}=\frac{F^u\Omega^v_\cH\cap
d^{-1}(F^{u+1}\Omega^{v+1}_\cH)}{F^{u+1}\Omega^v_\cH+d(F^u\Omega^{v-1}_\cH)}\;.$$
Since clearly $F^u\Omega^u_\cH=\Omega^u_{\cH,i=0}$,
we get $E_0^{\cdot,0}\equiv\Omega^\cdot_{\cH,i=0}$. 
Moreover,
\be\label{eqn:theta=0}\theta_X(\Omega^\cdot_{\cH,i=0})=0\ee
if the vector field $X$ is tangent to the $\cH$-orbit closures. Indeed it is easy to check
that Eq.~(\ref{eqn:theta=0}) follows if it is proved for $\hat{\cH}$. But since $\hat{\cH}$
preserves a parallelism, it is enough to prove Eq.~(\ref{eqn:theta=0}) for
$\hat{\cH}$-invariant functions, and Eq.~(\ref{eqn:theta=0}) is obvious in this case since 
such functions are constant on the $\hat{\cH}$-orbit closures.  

From Eq.~(\ref{eqn:theta=0}) we get $i_Xd(\Omega^\cdot_{\cH,i=0})=0$ for such $X$; i.e.
$d(\Omega^\cdot_{\cH,i=0})\subset\Omega^\cdot_{\cH,i=0}$, yielding
 \be\label{eqn:i=0}
E_1^{\cdot,0}\equiv\Omega_{\cH, i =0}\cong\Omega(W)_{{\rm O}(n), i =0}\;,\ee 
where the isomorphism is given by
$\alpha\mapsto\bar{\alpha}$ if $\pi^\ast\alpha=\pi_b^\ast\bar{\alpha}$. Therefore, from
Eq.~(\ref{eqn:W}) and the result in \cite{Verona} we get 
\be\label{eqn:E2} E_2^{\cdot,0}\cong
H^\cdot(T/\bar{\cH})\;.\ee 
Now, from the general theory of spectral sequences, there is a canonical injection
$E_2^{1,0}\ar H^1(T)_\cH$. So Eq.~(\ref{eqn:E2}) yields an injection 
$H^1(T/\bar{\cH})\ar H^1(T)_\cH$. 

If $F^kH^\cdot(T)_\cH$ denotes the filtration
of $H^\cdot(T)_\cH$ induced by the filtration of $\Omega^\cdot_\cH$, then
$E_\infty^{0,1}\equiv H^1(T)_\cH/F^1H^1(T)_\cH$ and $E_\infty^{1,0}\equiv F^1H^1(T)_\cH\cong
H^1(T/\bar{\cH})$ by Eq.~(\ref{eqn:E2}). Moreover
$\pi_\infty^\ast:E_\infty^{0,1}(\cH)\ar E_\infty^{0,1}(\hat{\cH})$ is  injective since
$\pi^\ast:H^1(T)_\cH\ar H^1(\hat{T})_{\hat{\cH}}$ is easily checked to be injective with
usual arguments involving the standard spectral sequence defined by $\pi$. 

A vector bundle $\rho:V\ar T$ will be called an $\cH$-vector bundle if, for any
diffeomorphism $h:U_1\ar U_2$ in $\cH$, there is a vector bundle homomorphism
$\tilde{h}:\rho^{-1}(U_1)\ar\rho^{-1}(U_2)$ over $h$ satisfying $\widetilde{{\rm
id}_T}={\rm id}_V$,  $\widetilde{h_1h_2}=\widetilde{h_1}\widetilde{h_2}$, and
$\tilde{h}|_{\rho^{-1}(U)}=\widetilde{h|_U}$ for every open subset $U\subset{\rm dom}\,h$. A
connection on $V$ will be called an $\cH$-connection if it is invariant by the $\tilde{h}$.
The following is a natural example of an $\cH$-vector bundle with an $\cH$-flat connection.
The locally trivial sheaf of infinitesimal transformations of $\bar{\cH}$ will be
denoted by $\cC=\cC(\cH)$ \cite{Salem}. Such $\cC$ is a sheaf of Lie algebras, whose typical
fiber is the opposite Lie algebra $\g^-$ of $\g$. The corresponding $\cH$-vector bundle will
be denoted by $C= C(\cH)$, and the corresponding $\cH$-flat connection by $\nabla$. By
naturality, the multiplicative homomorphism $\Delta_\ast:H^\cdot({\rm gl}(n),{\rm O}(n))\ar
H^\cdot(T)$ defined by $\nabla$ is a composite of a multiplicative homomorphism
$\Delta_\ast=\Delta(\cH)_\ast:H^\cdot({\rm gl}(n),{\rm O}(n))\ar H^\cdot(T)_\cH$ and the
canonical homomorphism $H^\cdot(T)_\cH\ar H^\cdot(T)$. This yields secondary  characteristic
invariants  $\Delta_\ast(y_i)=\Delta(\cH)_\ast(y_i)\in H^{2i-1}(T)_\cH$.  Indeed the
representatives defined in \cite{KamberTond1} are $\cH$-invariant; in particular 
$\Delta(y_1)\in \Omega^1_\cH$. $\Delta(y_1)$ and $\Delta_\ast(y_1)\in H^1(T)_\cH$ are the
objects of our study.  

The results in
\cite{Sergiescu} have obvious versions for complete pseudogroups of local isometries. In
particular, when $T$ has an $\cH$-invariant orientation, the top-dimensional invariant
cohomology $H^n(T)_\cH$ is non-trivial if and only if $\bigwedge^m\cal C$ is a trivial
sheaf, which is equivalent to  $\Delta_\ast(y_1)=0$. 

\section{The form $\Delta(y_1)\in\Omega^1(T)_\cH$ when $\cH$ preserves a parallelism}
\label{sec:Delta}
With the notation of Sect.~\ref{sec:preliminaries}, suppose $\cH$ preserves a parallelism on
$T$. Then, as a particular case of a foliation in a Riemannian manifold, there is a bigrading
of $\Omega$ given by the fibers of $\pi_b$: If $\cal V$ is the vertical bundle of $\pi_b$ and
$\cal Q$ the orthogonal complement of $\cal V$, then $$\Omega^{u,v}=\Cinf(\bigwedge^u{\cal
Q}^\ast\otimes\bigwedge^v{\cal V}^\ast)\; ,\quad u,v\in\Z\; .$$ Such bigrading is
$\cH$-invariant, and thus restricts to $\Omega_\cH$. The de~Rham derivative decomposes as
$d=d_{0,1}+d_{1,0}+d_{2,-1}$, where each $d_{i,j}$ is bihomogeneous of bidegree $(i,j)$, and
the usual formulae are satisfied (see e.g. \cite{Alv1}). Clearly
$F^k\Omega_\cH=\Omega_\cH^{k,\cdot}\wedge\Omega_\cH$, yielding canonical identities 
$(E_0,d_0)\equiv(\Omega_\cH,d_{0,1})$ and
$(E_1,d_1)\equiv(H((\Omega_\cH,d_{0,1}),d_{0,1\ast})$. 

Let
$\chi\in\Omega_\cH^{0,m}$ be the characteristic form of the fibers of $\pi_b$, where $m=\dim
\g$. There is a form $\tau\in\Omega_\cH^{0,1}$ such that, for any $\cH$-invariant vector
field $Y$ tangent to the $\cH$-orbit closures,    \be\label{eqn:tau}\theta_Y\chi\in
-\tau(Y)\,\chi+F^1\Omega_\cH^m\; .\ee Indeed, if $\cH_F$ is the restriction of $\cH$ to any
$\cH$-orbit closure $F$, then the Lie algebra ${\cal X}_F$ of $\cH_F$-invariant vector fields
on $F$ is isomorphic to $\g$ by Salem's description, and the restriction $\tau_F$ of $\tau$
to $F$ is the trace of the adjoint representation of ${\cal X}_F$. So $\tau=0$ if and only if
$\g$ is unimodular.  On the other hand, the mean curvature form $\kappa$ of the $\cH$-orbit
closures is in $\Omega_\cH^{1,0}=F^1\Omega_\cH^1$, and satisfies Rummler's formula 
\be\label{eqn:kappa} \theta_Z\chi\in -\kappa(Z)\,\chi+F^1\Omega_\cH^m\ee
for any $\cH$-invariant vector field $Z$ orthogonal to the orbit closures. Also, with the
notation of Sect.~\ref{sec:remark} for $V= C$, let $\kappa_C$ be the $\hat{g}$-mean
curvature of the fibers of the projection of $C$ to $T$.

\bp \label{prop:tau+kappa} $\tau+\kappa=2\pi\,\Delta(y_1)$. \ep

{\em Proof}.  For any fixed subset $U\subset T/\bar{\cH}$, let $ C_U$ be the restriction
of $ C$ to $\pi_b^{-1}(U)$. The statement of this result is a local property, thus it is
enough to prove it on $\pi_b^{-1}(U)$. Then, by Salem's description, we can assume
$\pi_b^{-1}(U)=U\times G$, where $\bar{\cH}$ is generated by the action of $G$, acting
by left multiplication on itself and trivially on $U$. Thus $ C_U\equiv U\times
G\times\g^-\equiv U\times{\rm T}G\subset{\rm T}(U\times G)$, and $\Cinf( C_U)\equiv
\Cinf(U\times G,\g^-)$, where the $\nabla$-parallel sections of $ C_U$ are identified with
the constant functions with values in $\g^-$. 

Let $X$ be a unitary section of $\bigwedge^m C_U$, which can be
considered as a function on $U\times G$ with values in $\bigwedge^m\g^-$. By
Eqs.~(\ref{eqn:tau}) and~(\ref{eqn:kappa}) we have 
\be\label{eqn:tau+kappa}\theta_YX=(\tau+\kappa)(Y)\,X\ee
for any $\cH$-invariant vector field $Y$ on $U\times G$.  

The vertical bundle of $ C_U$ can be canonically identified with the trivial bundle ${\bf
C}_U\times\g^-$. Hence, using the notation of the proof of Proposition~\ref{prop:kappa} for
$V= C_U$, $\hat{X}$ can be considered as a function on $ C_U$ with values in
$\bigwedge^m\g^-$, which is clearly equal to the pull-back of $X$. So, as in
Eq.~(\ref{eqn:Rummler}),  \be\label{eqn:hat} \widehat{\theta_YX}=\theta_{\tilde{Y}}\hat{X}=
\kappa_C(\tilde{Y})\,\hat{X}\ee
for any vector field $Y$ on $U\times G$, where $\tilde{Y}$ is the horizontal lift of $Y$. 
Then the result follows from Eqs.~(\ref{eqn:tau+kappa}) and (\ref{eqn:hat}), and
Proposition~\ref{prop:kappa}.  $\Box$  

\section{Unimodularity of the structural Lie algebra}\label{sec:unimodularity}
\bt\label{thm:unimodularity} If $\cH$ is a complete pseudogroup of local isometries and $\g$
its structural Lie algebra, then $\bar{\Delta}_\ast(y_1)=0$ if and only if $\g$ is
unimodular.\et

{\em Proof}.  Since $\pi_\infty^\ast:E_\infty^{0,1}(\cH)\ar E_\infty^{0,1}(\hat{\cH})$ is
injective and  $\pi_\infty^\ast\bar{\Delta}(\cH)_\ast(y_1)=
\bar{\Delta}(\hat{\cH})_\ast(y_1)$, we can assume $\cH$ preserves a parallelism. 

With the above assumption, if $\bar{\Delta}_\ast(y_1)=0$ then $\Delta_\ast(y_1)\in
F^1H^1(T)_\cH$. So there is some $\cH$-invariant function $f$ such that $\tau+\kappa+df\in
F^1\Omega_\cH^1$. But $\kappa+df\in F^1\Omega_\cH^1$ because $d(F^1\Omega_\cH^1)=0$.
Therefore $\tau\in\Omega_\cH^{0,1}\cap F^1\Omega_\cH^1=0$, and $\g$ is unimodular. 

If $\g$ is unimodular, then $\tau=0$ and thus $\kappa\in F^1\Omega_\cH^1$ represents
$2\pi\,\Delta_\ast(y_1)$. So $\Delta_\ast(y_1)\in F^1H^1(T)_\cH$ and 
$\bar{\Delta}_\ast(y_1)=0$. $\Box$  

From Theorem~\ref{thm:unimodularity}, any vanishing result for $H^1(T)_\cH$ yields the
unimodularity of $\g$. Such a result is proved e.g. in \cite{Alv8} by using Morse
inequalities for pseudogroups of local isometries.

\section{Minimality of the orbit closures when $\cH$ preserves a parallelism}
\label{sec:minimality}
Assume $\cH$ preserves a parallelism on $T$. With the notation of Sect.~\ref{sec:Delta}, let
$\nu$ denote the normal bundle of the fibers of $\pi_b$, which is canonically isomorphic
to $\cal Q$, and let $\Cinf(\nu^\ast\otimes{\cal V})_\cH$ be the space of $\cH$-invariant
sections of the $\cH$-vector bundle $\nu^\ast\otimes{\cal V}$. For
$\sigma\in\Cinf(\nu^\ast\otimes{\cal V})_\cH$ and $x\in T$, let ${\cal
Q}_x^\sigma=\{v+\sigma_x(\bar{v}):\ v\in{\cal Q}_x\}$, where $\bar{v}$ is the
element in $\nu$ defined by $v$. Clearly ${\cal Q}^\sigma=\bigsqcup_{x\in T}{\cal
Q}^\sigma_x$ is an $\cH$-subbundle of the tangent bundle of $T$ which is complementary of
$\cal V$. The correspondence $\sigma\mapsto{\cal Q}^\sigma$ defines a bijection between 
$\Cinf(\nu^\ast\otimes{\cal V})_\cH$ and the $\cH$-bundles of tangent vectors which
are complementary of $\cal V$. This correspondence depends on $\cal Q$, and thus on the given
metric $g$. For each such $\sigma$, there is a unique $\cH$-invariant metric $g^\sigma$ on
$T$ such that $g$ and $g^\sigma$ induce the same metric on $T/\bar{\cH}$ and  define the
same metric on $\cal V$, and so that the $g^\sigma$-orthogonal complement of $\cal V$ is
${\cal Q}^\sigma$. The $g^\sigma$-mean curvature form of the orbit closures will be
denoted by $\kappa^\sigma$. 

For any $\cH$-invariant function $h$ on $T$, consider also the orthogonal sum $g^h$ of the
restriction of $g$ to $\cal Q$ and the restriction of ${\rm e}^h\,g$ to $\cal V$. Such $g^h$
will be said to be obtained from $g$ by a scalar change along the orbit closures.

\bt\label{thm:minimality} Let $\cH$ be a complete pseudogroup of local isometries which
preserves a parallelism on a manifold $T$, and $\g$ the structural Lie algebra. Then:
\begin{itemize}
\item[(i)] If $\g$ is unimodular, $2\pi\,\Delta_\ast(y_1)$ is the class of all possible mean
curvature forms of the orbit closures for all the $\cH$-invariant metrics on $T$. Thus, in
this case, $\Delta_\ast(y_1)=0$ if and only if the orbit closures are
minimal submanifolds for some $\cH$-invariant metric.
\item[(ii)] If $\g$ is not unimodular, any element in $\Omega_\cH^{1,0}$ is the mean curvature
form of the orbit closures for some $\cH$-invariant metric on $T$. Thus, in this case, the
orbit closures are minimal submanifolds for some $\cH$-invariant metric. 
\end{itemize}
Thus $\Delta_\ast(y_1)=0$ if and only if $\g$ is unimodular and the orbit closures
are minimal submanifolds for some $\cH$-invariant metric on $T$.\et
 
{\em Proof}.  Let $g$ be any
$\cH$-invariant metric. If $\g$ is unimodular, $\tau=0$ and $\kappa$ represents
$2\pi\,\Delta_\ast(y_1)$. On the other hand, any element in this class can be realized as
the mean curvature form of the orbit closures for some metric obtained from $g$ by a scalar
change along the orbit closures. 

Suppose $\g$ is not unimodular, and thus $\tau$ is a non-vanishing form. So there is an
$\cH$-invariant $\pi_b$-vertical vector field $Y$ such that $\tau(Y)=1$. ($Y$ can be chosen
to be $|\tau|^{-2}$ times the $g$-dual vector field of $\tau$.) Take any
$\alpha\in\Omega_\cH^{1,0}$ and define $\sigma\in\Cinf(\nu^\ast\otimes{\cal V})_\cH$ by
$\sigma_x(\bar{v})=(\alpha+\kappa)_x(v)\,Y_x$ for any $x\in T$ and any tangent vector
$v$ at $x$, and where $\bar{v}$ is the element defined by $v$ in $\nu_x$.  Let
$P^\sigma$ denote the $g^\sigma$-orthogonal projection of the tangent bundle of $T$ onto $\cal
V$. It is easily verified that
\be\label{eqn:P} P^\sigma(v)=(\alpha+\kappa)_x(v)\,Y_x\quad\mbox{for}\quad v\in{\cal
Q}_x\; .\ee

For a local orientation of the fibers of $\pi_b$, let $\chi$ and $\chi^\sigma$ be the
corresponding characteristic forms for $g$ and $g^\sigma$, respectively. Let
$Y_1,\ldots,Y_m$ be a local orthonormal frame of $\cal V$ such that $Y_1=Y/|Y|$. By
Eq.~(\ref{eqn:P}), for any $X\in\Cinf({\cal Q})$ we have
\ba \chi^\sigma(X\wedge Y_2\wedge\ldots\wedge Y_m)&=&\chi\left( P^\sigma(X)\wedge
Y_2\wedge\ldots\wedge Y_m\right)\\
&=&\chi\left( (\alpha+\kappa)(X)\,Y\wedge Y_2\wedge\ldots\wedge Y_m\right)\\
&=&(\alpha+\kappa)(X)\,|Y|\; ,\ea
$$\chi^\sigma(Y_1\wedge\ldots\wedge Y_{i-1}\wedge X\wedge Y_{i+1}\wedge\ldots\wedge
Y_m)=0\; ,$$ 
$$\chi^\sigma(Y_1\wedge\ldots\wedge Y_m)=1\; .$$
Thus
\be\label{eqn:chi}\chi^\sigma=\chi+(\alpha+\kappa)\wedge i_Y\chi\; .\ee
Therefore
\ba d\chi^\sigma&=&d\left(\chi+\alpha+\kappa)\wedge i_Y\chi\right)\\
&=&d\chi+d(\alpha+\kappa)\wedge i_Y\chi-(\alpha-\kappa)\wedge d\,i_Y\chi\\
&\in&-\kappa\wedge\chi+d(\alpha+\kappa)\wedge i_Y\chi-
(\alpha+\kappa)\wedge(\theta_Y- i_Y d)\chi+F^2\Omega_\cH^m\\
&=&-\kappa\wedge\chi-(\alpha+\kappa)\wedge\theta_Y\chi+F^2\Omega_\cH^m\; .
\ea
But $\theta_Y\chi=-\chi$ by the definition of $\tau$ and the choice of $Y$. So
$$d\chi^\sigma\in -\alpha\wedge\chi+F^2\Omega_\cH=
-\alpha\wedge\chi^\sigma+F^2\Omega_\cH\; ,$$ yielding $\kappa^\sigma=\alpha$ by Rummler's mean
curvature formula, and the proof is finished.  $\Box$

\begin{example}[Y.~Carri\`ere]\label{ex:Carriere} The affine Lie group $A$ can be identified
with $\R^2$ with the group structure given by $(t,s)(t',s')=(t+t',\lambda^ts'+s)$ for any
fixed $\lambda>1$. Consider the pseudogroup generated by the left action of the closed
subgroup $K=\Z\times\R\subset A$ on $A$. Clearly $K\backslash A\equiv{\Bbb S}^1$, the
structural Lie algebra is abelian, and we have  $H^2(A)_K=0$
\cite{Carriere}. Hence $\bar{\Delta}_\ast(y_1)=0$ and $\Delta_\ast(y_1)\neq 0$.
Therefore there is no $K$-invariant metric on $A$ such that the right translates of $K$ are
minimal submanifolds. \end{example}

E.~Mac\'{\i}as and E.~Sanmart\'{\i}n have proved the
following \cite{MaciasSanmartin}: {\em If $H$ is a Lie subgroup of a Lie group $G$, and $H_0$
the connected component of $H$ which contains the identity element, then the right translates
of $H$ are minimal submanifolds for some metric on $G$}. Moreover, from the proof in
\cite{MaciasSanmartin} it can be easily seen that the above metric can be chosen to be
invariant by the left action of $H_0$. So the non-triviality
of $\Delta_\ast(y_1)$ in Example~\ref{ex:Carriere} depends on the disconnectedness of $K$. 
From Theorems~\ref{thm:unimodularity}
and~\ref{thm:minimality}, we get the following
generalization of the results in \cite{MaciasSanmartin}, where $\Delta_\ast(y_1)$ and
$\bar{\Delta}_\ast(y_1)$ are defined by the pseudogroup generated by the left action of
$H$ on $G$.

\bc With the above notation, we have:
\begin{itemize}
\item[(i)] $\bar{\Delta}_\ast(y_1)$ if and only if $H$ is unimodular.
\item[(ii)] If $H$ is unimodular, $\Delta_\ast(y_1)=0$ if and only if the right translates of
$H$ are minimal submanifolds for some metric on $G$ which is invariant by the left
$H$-action. 
\item[(iii)] If $H$ is not unimodular, the right translates of $H$ are minimal
submanifolds for some metric on $G$ which is invariant by the left $H$-action.
\end{itemize}\ec

If $\cH$ is not required to preserve a parallelism on $T$, then
Theorems~\ref{thm:unimodularity} and~\ref{thm:minimality} have the following consequence by
considering $\hat{\cH}$.

\bc\label{cor:general} If $\cH$ is a complete pseudogroup of local isometries, and $\g$ its
structural Lie algebra, then:
\begin{itemize}
\item[(i)] If $\g$ is unimodular, $\Delta_\ast(y_1)=0$ if and only if the $\hat{\cH}$-orbit
closures are minimal submanifolds for some $\hat{\cH}$-invariant metric.
\item[(ii)] If $\g$ is not unimodular, the $\hat{\cH}$-orbit
closures are minimal submanifolds for some $\hat{\cH}$-invariant metric.
\end{itemize}
Thus $\Delta_\ast(y_1)=0$ if and only if $\g$ is unimodular and the $\hat{\cH}$-orbit
closures are minimal submanifolds for some $\hat{\cH}$-invariant metric.\ec

\section{Application to Riemannian foliations}
Let $\F$ and $M$ be as in Sect.~\ref{sec:intro}. Let $\cH$ be the representative of the
holonomy pseudogroup of $\F$ canonically defined on a manifold $T$ by some regular covering
of $M$ (see e.g. \cite{Haefliger80,Haefliger85}). Then any fixed transverse
Riemannian structure of $\F$ canonically corresponds to an $\cH$-invariant metric on $T$ so
that $\cH$ is a complete pseudogroup of local isometries, and there are  canonical
isomorphisms    \be\label{iso1} \Omega^\cdot(T)_\cH\cong\Omega^\cdot(M/\F)\; ,\ee 
\be\label{iso2} H^\cdot(T)_\cH\cong H^\cdot(M/\F)\; .\ee

More precisely, let $\{(U_i,f_i)\}$ be a regular covering of $M$. The restriction of $\F$
to each $U_i$ is given by the submersion $f_i$ of $U_i$ onto some manifold $T_i$. The
regularity of this covering means that there are well defined diffeomorphisms
$h_{ij}:f_i(U_i\cap U_j)\ar f_j(U_i\cap U_j)$ such that $h_{ij}f_i=f_j$ on $U_i\cap
U_j$. Then $T=\bigsqcup_iT_i$ and $\cH$ is generated by the $h_{ij}$. The metric on $T$ is
determined by requiring the $f_i$ to be Riemannian submersions. The isomorphism in 
Eq.~(\ref{iso1}) is given by $\alpha\mapsto \alpha'$ where
$f_i^\ast(\alpha|_{T_i})=\alpha'|_{U_i}$. Moreover each $f_i^\ast C(\cH)$ is canonically
isomorphic to the restriction of $C(\F)$ to $U_i$, so the classes $\Delta(\cH)_\ast(y_1)$ and
$\Delta(\F)_\ast(y_1)$ correspond to each other by Eq.~(\ref{iso2}). 

{\em Proof of Theorem~\ref{thm:xi}}. We can suppose $\F$ is transversely
parallelizable by using $\hat{\F}$ in a standard way. Then
$M/\bar{\F}\equiv T/\bar{\cH}$ is a manifold, and the canonical map $\pi_b:M\ar
M/\bar{\F}$ a fiber bundle projection whose fibers define thus a foliation
$\bar{\F}$. We can suppose the metric on $M$ is chosen so that the leaf closures are
minimal submanifolds. On any fixed $U_i$, the vector bundles ${\rm T}\F$,
${\rm T}\bar{\F}$ and ${\rm T}\F^\perp\cap{\rm T}\bar{\F}$ are orientable,
and take the unitary sections $X$, $X'$ and $X''$ on $U_i$ defining respective positive
orientations of $\bigwedge^p{\rm T}\F$, $\bigwedge^{p+m}{\rm T}\bar{\F}$ and 
$\bigwedge^m\left({\rm T}\F^\perp\cap{\rm T}\bar{\F}\right)$, where $p=\dim\F$  
(thus $p+m=\dim\bar{\F}$). The orientations can be chosen so that $X'=X\wedge X''$. Let
$Y$ be an infinitesimal transformation of $\F$ which is orthogonal to the leaves.
Since $X\wedge\bigwedge^+\Cinf({\rm T}\F)=0$, we get
$X\wedge\theta_YX''=(2\pi\,f_i^\ast\Delta(\F)(y_1))(Y)\,X'$ by Eq.~(\ref{eqn:tau+kappa}) and
Proposition~\ref{prop:tau+kappa}. Hence 
$$\theta_YX'=\left(\kappa_\F+2\pi\,
f_i^\ast\Delta(\cH)(y_1)\right)(Y)\,X'$$ 
by Rummler's formula, where $\kappa_\F$ is the mean curvature form of the leaves. On the one
hand, if $Y$ is orthogonal to $\bar{\F}$ , $\theta_YX'=0$ because the leaves of
$\bar{\F}$ are minimal submanifolds, and we get $(\kappa_\F+2\pi\,\Delta(\F)(y_1))(Y)=0$
on $M$. On the other hand,  if $Y$ is tangent to $\bar{\F}$,  $\theta_YX'=-{\rm
div}_b(Y)\,X'$  where ${\rm div}_b$ denotes the divergence on the fibers of $\pi_b$, and
we get $(\kappa_\F+2\pi\,\Delta(\F)(y_1))(Y)=-{\rm div}_b(Y)$ on $M$. Therefore the function
$(\kappa_\F+2\pi\,\Delta(\F)(y_1))(Y)$ has trivial integral on the fibers of $\pi_b$  
for any infinitesimal transformation $Y$ of $\F$ on $M$. This implies that 
$\kappa_\F+2\pi\,\Delta(\F)(y_1)$ has trivial basic component \cite{Alv7}; i.e. the basic
component of $\kappa_\F$ is equal to $-2\pi\,\Delta(\F)(y_1)$, and
the result follows. $\Box$   

The other results in the Sect.~\ref{sec:intro} now follow directly from 
Theorem~\ref{thm:xi} and the results in Sects.~\ref{sec:unimodularity}
and~\ref{sec:minimality}.

\end{document}